\documentclass[11pt]{article}

\usepackage{makeidx}
\usepackage{amsmath,amssymb,amsthm}

\newcommand{\lam}{\lambda}
\newcommand{\eset}{\emptyset}

\theoremstyle{definition}
\newtheorem{definition}{Definition}
\newtheorem{theorem}{Theorem}
\newtheorem{prop}[theorem]{Proposition}

\newtheorem{cor}[theorem]{Corollary}

\newtheorem{example}{Example}

\newcommand{\content}{\mathrm{content}}
\newcommand{\crr}{\mathrm{cr}}
\newcommand{\nee}{\mathrm{ne}}
\newcommand{\st}{\,|\,}
\newcommand{\zz}{\mathbb{Z}}
\newcommand{\rr}{\mathbb{R}}
\newcommand{\cc}{\mathbb{C}}
\newcommand{\sn}{\mathfrak{S}_n}
\newcommand{\sm}{\mathfrak{S}_m}

\title{Crossings and Nestings of  Matchings and Partitions}
\author{William Y.C. Chen$^{1,2}$, Eva Y.P. Deng$^3$, Rosena R.X. Du$^4$ \\
Richard P.  Stanley$^{5,6}$, and Catherine H. Yan$^{7,8}$
\vspace{.3cm} \\
$^{1}$Center for Combinatorics, LPMC \\
Nankai University, Tianjin 300071, P.R. China
 \vspace{.3cm} \\
$^3$Department of Applied Mathematics \\
Dalian University  of Technology, Dalian, Liaoning 116024, P.R.China
\vspace{.3cm} \\
$^4$Department of Mathematics \\ 
East China Normal University, Shanghai 200062, P.R. China
\vspace{.3cm} \\
$^5$Department of Mathematics \\
Massachusetts Institute of Technology, Cambridge, MA 02139
\vspace{.3cm} \\
$^7$Department of Mathematics\\
Texas A\&M University, College Station, TX 77843
\vspace{.3cm} \\
$^1$chen@nankai.edu.cn,  $^3$dengyp@eyou.com, $^4$rxdu@math.ecnu.edu.cn,
\\ $^5$rstan@math.mit.edu, $^7$cyan@math.tamu.edu
}

\begin{document}

\maketitle

\begin{abstract}
We present results on the enumeration of crossings and nestings for
matchings and set partitions.  Using a bijection between partitions and
vacillating tableaux, we show that if we fix the sets of minimal block
elements and maximal block elements, the crossing number and the
nesting number of partitions have a symmetric joint distribution.  It
follows that the crossing numbers and the nesting numbers are distributed
symmetrically over all partitions of $[n]$, as well as over all
matchings on $[2n]$.  As a corollary, the number of $k$-noncrossing
partitions is equal to the number of $k$-nonnesting partitions. The
same is also true for matchings. An application is given to the
enumeration of matchings with no $k$-crossing (or with no
$k$-nesting).
\end{abstract}



\footnotetext[2]{partially supported by the 973 Project on
Mathematical Mechanization, the National Science Foundation, 
the Ministry of Education  and
the Ministry of Science and Technology of China.}
\footnotetext[6]{Partially supported by NSF grant \#DMS-9988459.}

\footnotetext[8]{Partially supported by NSF grant \#DMS-0245526 and a
Sloan Fellowship.}

%
%
\section{Introduction}
A (complete) matching on $[2n]=\{1,2,\dots,2n\}$ is  a partition  of
$[2n]$ of type $(2,2, \dots, 2)$. It can be represented by listing its
$n$ blocks, as $\{(i_1, j_1), (i_2, j_2), \dots, (i_n, j_n)\}$
where $i_r < j_r$ for $ 1 \leq r \leq n$. Two blocks (also called
arcs)
$(i_r, j_r)$ and $(i_s, j_s)$ form a  \emph{crossing} if $i_r < i_s <
j_r < j_s$; 
they form a \emph{nesting} if $i_r < i_s < j_s < j_r$.
It is well-known that the number of matchings on $[2n]$ with no
crossings (or with no nestings) is given by the $n$-th Catalan
number
\[
C_n= \frac{1}{n+1} \binom{2n}{n}.
\]
See \cite[Exercise 6.19]{Stanley99} for many combinatorial interpretations
of Catalan numbers, where item (o) is for noncrossing matchings, and item
(ww) can be viewed as  nonnesting matchings, in which the blocks of
the matching are the  columns of the standard Young tableaux of shape
$(n,n)$. Nonnesting matchings are also one of the items of
\cite{Stanley_web}.

Let $k \geq 2$ be an integer.  A \emph{$k$-crossing} of a matching $M$ is a
set of $k$ arcs $(i_{r_1}, j_{r_1})$, $(i_{r_2}, j_{r_2}), \dots$,
$(i_{r_k}, j_{r_k})$ of $M$ such that $i_{r_1} <i_{r_2} < \cdots <
i_{r_k} < j_{r_1} < j_{r_2} < \cdots < j_{r_k}$.
 A matching without any
$k$-crossing is a \emph{$k$-noncrossing matching.}  Similarly, a
$k$-nesting is a set of $k$ arcs $(i_{r_1}, j_{r_1}), (i_{r_2},
j_{r_2})$, $\dots$, $(i_{r_k}, j_{r_k})$ of $M$ such that $i_{r_1}
<i_{r_2}$ $<$ $\cdots$ $<$ $i_{r_k} < j_{r_k}< \cdots < j_{r_2} < j_{r_1}$.
A matching without any $k$-nesting is a \emph{$k$-nonnesting matching.}

Enumeration on crossings/nestings of matchings has been studied
for the cases $k=2$ and $k=3$. For  $k=2$, in addition to the
above results on Catalan numbers, the distribution of the number
of $2$-crossings has been studied by Touchard \cite{Touchard52},
and later more explicitly by Riordan \cite{Riordan75},  who gave a
generating function. M. de Sainte-Catherine \cite{SC83} proved
that $2$-crossings and $2$-nestings are  identically distributed
over all matchings of $[2n]$, i.e., the number of matchings with
$r$ $2$-crossings is equal to the number of matchings with $r$
$2$-nestings.

The enumeration of 3-nonnesting matchings was first studied by
Gouyou-Beauschamps \cite{GB89}, in which he gave a bijection
between involutions with no decreasing sequence of length 6 and pairs
of noncrossing Dyck left factors by a recursive construction. His
bijection is essentially a correspondence between 3-nonnesting
matchings and pairs of noncrossing Dyck paths, where a matching can
also be considered as a fixed-point-free involution. We observed that
the number of 3-noncrossing matchings also equals the number of
pairs of noncrossing Dyck paths, and a one-to-one correspondence
between 3-noncrossing matchings and pairs of noncrossing Dyck paths
can be built recursively.

In this paper, we extend the above results. Let $\crr(M)$ be
maximal $i$ such that $M$ has an $i$-crossing, and $\nee(M)$ the
maximal  $j$ such that $M$ has a $j$-nesting. Denoted by
$f_n(i,j)$ the number of matchings $M$ on $[2n]$ with $\crr(M)=i$
and $\nee(M)=j$. We shall prove that $f_n(i,j)=f_n(j,i)$. As a
corollary, the number of matchings on $[2n]$ with $\crr(M)=k$
equals the number of matchings $M$ on $[2n]$ with $\nee(M)=k$.

Our construction applies to a more general structure, viz., partitions
of a set. Given a partition $P$ of $[n]$, denoted by $P \in \Pi_n$, we
represent $P$ by a graph on the vertex set $[n]$ whose edge set
consists of arcs connecting the elements of each block in numerical
order.  Such an edge set is called the \emph{standard representation}
of the partition $P$.  For example, the standard representation of
1457-26-3 is $\{(1,4), (4, 5), (5, 7), (2,6)\}$.  Here we always write
an arc $e$ as a pair $(i,j)$ with $i < j$, and say that $i$ is the
\emph{lefthand endpoint} of $e$ and $j$ is the \emph{righthand
endpoint} of $e$.

Let $k \geq 2$ and $P \in \Pi_n$.  Define a \emph{$k$-crossing} of
$P$ as a $k$-subset $(i_1, j_1), (i_2, j_2), \dots, (i_k, j_k)$ of
the arcs in the standard representation of $P$ such that $i_1 <
i_2 < \cdots < i_k < j_1 < j_2 < \cdots < j_k$.  Let $\crr(P)$ be
the maximal $k$ such that $P$ has a $k$-crossing.  Similarly,
define a \emph{$k$-nesting} of $P$ as a $k$-subset $(i_1, j_1),
(i_2, j_2), \dots, (i_k, j_k)$ of the set of arcs in the standard
representation of $P$ such that $i_1 < i_2 < \cdots < i_k < j_k <
\cdots < j_2 < j_1$, and $\nee(P)$ the maximal $j$ such that $P$
has a $j$-nesting. Note that when restricted to complete
matchings, these definitions agree with the ones given before.

Let $g_n(i,j)$ be the number of partitions  $P$ of $[n]$ with
$\crr(P)=i$ and $\nee(P)=j$. We shall prove that $g_n(i,j)=g_n(j,i)$, for
all $i,j$ and $n$. In fact, our result is much stronger. We present
a generalization  which
implies the symmetric distribution
of $\crr(P)$ and $\nee(P)$ over all partitions in $\Pi_n$,
as well over all complete matchings on $[2n]$.

To state the main result, we need some notation.
Given $P \in \Pi_n$, define
\begin{eqnarray*}
\min(P)=\{ \text{minimal block elements of $P$}\}, \\
\max(P)=\{ \text{maximal block elements of $P$}\}.
\end{eqnarray*}
For example, for $P=\mbox{135-26-4}$, $\min(P) =\{1,2,4\}$ and
$\max(P)=\{4,5,6\}$. The pair $(\min(P), \max(P))$ encodes some useful
information about the partition $P$. For example, the number of blocks
of $P$ is $|\min(P)|=|\max(P)|$; number of singleton blocks is
$|\min(P) \cap \max(P)|$; $P$ is a (partial) matching if and only if
$\min(P) \cup \max(P)=[n]$, and $P$ is a
complete matching if in addition, $\min(P) \cap \max(P)=\eset$.

Fix $S, T \subseteq [n]$ with $|S|=|T|$.  Let $P_n(S, T)$ be the
set $\{ P \in \Pi_n: \min(P)=S, \max(P)=T\}$, and
$f_{n, S, T}(i, j)$ be the cardinality of the set
$\{ P \in P_n(S, T): \crr(P)=i, \nee(P)=j\}$.

\begin{theorem} \label{thm1}
\begin{eqnarray}\label{eq1}
f_{n, S, T}(i,j)=f_{n, S, T}(j, i).
\end{eqnarray}
In other words,
\begin{eqnarray} \label{dist}
\sum_{P \in P_n(S, T)} x^{\crr(P)}y^{\nee(P)} =
\sum_{P \in P_n(S, T)} x^{\nee(P)}y^{\crr(P)}.
\end{eqnarray}
 That is, the statistics $\crr(P)$ and $\nee(P)$ have a symmetric
joint distribution over each set $P_n(S, T)$. $\Box$ 
\end{theorem}

Summing over all pairs $(S, T)$ in \eqref{eq1}, we get
\begin{eqnarray}\label{eq2}
g_n(i,j)=g_n(j,i). 
\end{eqnarray}
We say that a partition $P$ is \emph{$k$-noncrossing} if $\crr(P)<k$.
It is \emph{ $k$-nonnesting} if $\nee(P) <k$.
Let $\mathrm{NCN}_{k,l}(n)$ be the number of partitions of $[n]$ that are
$k$-noncrossing and $l$-nonnesting. Summing over $1 \leq i < k$ and
$1 \leq j <l $ in \eqref{eq2}, we get the following corollary.
\begin{cor}\label{ncn}
$\mathrm{NCN}_{k,l}(n)=\mathrm{NCN}_{l,k}(n). \qquad \Box$
\end{cor}

Letting $l>n $,  Corollary \ref{ncn} becomes the following result.
\begin{cor} \label{nc_nn}
$\mathrm{NC}_k(n)=\mathrm{NN}_k(n)$, where $\mathrm{NC}_k(n)$ is
the number of $k$-noncrossing partitions of $[n]$, and
$\mathrm{NN}_k(n)$ is the number of
$k$-nonnesting partitions of $[n]$. ~~~~$\Box$ 
\end{cor}

Theorem \ref{thm1} also applies to complete matchings.
A partition $P$ of $[2n]$
is a complete matching if and only  if $|\min(P)|=|\max(P)|=n$ and
$\min(P) \cap \max(P)=\eset$. (It follows that $\min(P) \cup \max(P)=[2n]$.)
Restricting Theorem \ref{thm1} to disjoint pairs  $(S, T)$ of $[2n]$ with
$|S|=|T|=n$, we get the following
result on the crossing and nesting number of complete matchings.

\begin{cor} \label{matching}
Let $M$ be a matching on $[2n]$. \\
1. The statistics $\crr(M)$ and $\nee(M)$ have a symmetric joint
distribution over $P_{2n}(S, T)$, where $|S|=|T|=n$, and $S, T$
are disjoint. \\
2. $f_n(i,j)=f_n(j,i)$ where $f_n(i,j)$ is the number of matchings on
$[2n]$ with $\crr(M)=i$ and $\nee(M)=j$. \\
3. The number of matchings on $[2n]$ that are  $k$-noncrossing and
$l$-nonnesting is equal to  the number of matchings on $[2n]$ that are
$l$-noncrossing and $k$-nonnesting.   \\
4. The number of
$k$-noncrossing matchings on $[2n]$ is equal to the number of
$k$-nonnesting matchings on $[2n]$.

\hfill $\Box$
\end{cor}

The paper is  arranged as follows.
In Section 2 we introduce the concept of vacillating tableau of
general shape, and give a bijective proof for the number of
vacillating tableaux of shape $\lam$ and length $2n$.
In Section 3 we apply the bijection of Section 2 to vacillating
tableaux of empty shape, and characterize crossings and nestings of
a partition by the corresponding vacillating tableau.
 The involution on the
set of vacillating tableaux defined by taking the conjugate to each
shape leads to an involution on partitions which
exchanges the statistics $\crr(P)$ and $\nee(P)$ while preserves
$\min(P)$ and $\max(P)$,  thus proving
Theorem \ref{thm1}.
Then we modify the bijection between partitions and vacillating
tableaux by taking isolated points into consideration, and give
an analogous result on the enhanced crossing number and nesting
number. This is the content of Section 4.
Finally in Section 5 we restrict our bijection to the set of complete
matchings and oscillating tableaux, and
study the enumeration of $k$-noncrossing matchings.
In particular, we construct bijections from $k$-noncrossing matchings for $k=2$
or $3$  to Dyck
paths  and pairs of noncrossing Dyck paths, respectively, and
present the generating function for the number of
$k$-noncrossing matchings.

%
%
\section{A Bijection between Set Partitions and Vacillating Tableaux}
\label{sec2}
Let $Y$ be \emph{Young's lattice}, that is, the set of all partitions
of all integers $n \in \mathbb{N}$ ordered component-wise, i.e.,
$(\mu_1, \mu_2, \dots) \leq (\lam_1, \lam_2, \dots, )$ if $\mu_i \leq
\lam_i$ for all $i$.  We write $\lam\vdash k$ or $|\lam|=k$ if $\sum
\lambda_i=k$. A vacillating tableau is a walk on the Hasse diagram of
Young's lattice subject to certain conditions.  The main tool in our
proof of Theorem~\ref{thm1} is a bijection between the set of set
partitions and the set of vacillating tableaux of empty shape
$\emptyset$.

\begin{definition}
A \emph{vacillating tableau} $V_\lam^{2n}$ of shape $\lam$ and
length $2n$ is a sequence $\lam^0, \lam^1, \dots, \lam^{2n}$ of integer
partitions such that  (i) $\lam^0=\eset$, and
$\lam^{2n}=\lam$, (ii) $\lam^{2i+1}$ is obtained from $\lam^{2i}$
by doing nothing (i.e., $\lam^{2i+1}=\lam^{2i}$) or deleting a
square, and (iii) $\lam^{2i}$ is obtained from $\lam^{2i-1}$ by
doing nothing or adding a square.
\end{definition}
In other words, a vacillating tableau of shape $\lam$ is a walk
on the Hasse diagram of Young's lattice from $\eset$ to $\lam$
where each step consists of either
(i) doing nothing twice,
(ii) do nothing then adding a square,
(iii) removing  a square then doing nothing, or
(iv) removing a square and then adding a square.
Note that if the length is larger than $0$,
$\lam^1=\eset$. If the vacillating tableau is of empty shape,
then $\lam^{2n-1}=\eset$ as well.

\begin{example} \label{VT1}
Abbreviate $\lam=(\lam_1, \lam_2, \dots)$ by $\lam_1\lam_2\cdots$.
There are 5 vacillating  tableaux of shape $\emptyset$ and length
$6$. They are
\begin{eqnarray} \label{VT1-1}
\begin{array}{ccccccc}
\lam^0     & \lam^1& \lam^2   & \lam^3  & \lam^4&\lam^5 & \lam^6 \\
\emptyset  & \eset &  \eset   & \eset   & \eset & \eset & \eset  \\
\eset      & \eset &  \eset   & \eset   & 1     & \eset &\eset  \\
\eset      & \eset &  1       & \eset   & 1     & \eset & \eset \\
\eset      & \eset &  1       & \eset   & \eset     & \eset &\eset  \\
\eset      & \eset & 1        & 1       & 1  & \eset  &  \eset
\end{array}
\end{eqnarray}
\end{example}

\begin{example}
An example of a vacillating tableau of shape $11$ and length $10$ is
given by
\[
\eset, \eset, 1, 1, 1, 1, 2, 2, 21, 11, 11.
\]
\end{example}

\begin{theorem} \label{VT-number}
(i) Let $g_\lam(n)$ be the number of vacillating tableaux of shape
$\lam\vdash k$ and  length $2n$. By a \emph{standard Young
tableau} (SYT) of shape $\lambda$, we mean an array $T$ of shape
$\lambda$ whose entries are distinct positive integers that
increase in every row and column. The \emph{content} of $T$ is the
set of positive integers that appear in it. (We don't require that
content$(\lambda)=[k]$, where $\lambda\vdash k$.) We then have
\[
g_\lam(n) = B(n, k) f^\lam,
\]
where $f^\lam$ is the number of SYT's of shape $\lam$ and content
$[k]$, and $B(n, k)$ is the number of partitions of $[n]$ with $k$
blocks distinguished.

(ii) The exponential generating function of $B(n, k)$ is given by
\begin{eqnarray}
\sum_{n \geq 0} B(n, k) \frac{x^n}{n!} =
\frac{1}{k!}(e^x-1)^k \exp(e^x-1).
\end{eqnarray}
\end{theorem}

To prove Theorem \ref{VT-number}, we construct a bijection between the
set ${\cal V}_\lam^{2n}$ of vacillating tableaux of shape $\lam$ and
length $2n$, and pairs $(P, T)$, where $P$ is a partition of $[n]$,
and $T$ is an SYT of shape $\lam$ such that $\content(T) \subseteq
\max(P)$.  In the next section we apply this bijection to vacillating
tableaux of empty shape, and related it to the enumeration of crossing
and nesting numbers of a partition.  In the following we shall assume
familiarity with the RSK algorithm, and use row-insertion $P
\longleftarrow k$ as the basic operation of the RSK algorithm.  For
the notation, as well as some basic properties of the RSK algorithm,
see e.g.\ \cite[Chapter 7]{Stanley99}.  In general we shall apply the
RSK algorithm to a sequence $w$ of distinct integers, denoted by $w
\stackrel{\text{RSK}}{\longmapsto}(A(w), B(w))$, where $A(w)$ is the
(row)-insertion tableau and $B(w)$ the recording tableau.  The shape
of the SYT's $A(w)$ and $B(w)$ is also called the \emph{shape} of the
sequence $w$.

{\bf The Bijection $\psi$ from Vacillating Tableaux to Pairs $(P, T)$}.  \\
Given a vacillating tableau $V=(\eset=\lam^0, \lam^1, \dots,
\lam^{2n}=\lam)$, we will recursively define a sequence $(P_0,
T_0)$, $(P_1, T_1)$,$\dots$, $(P_{2n}, T_{2n})$ where $P_i$ is a
set of ordered pairs of integers in $[n]$, and $T_i$ is an SYT of shape
$\lam^i$.  Let $P_0$ be the empty set, and let $T_0$ be the empty
SYT (on the empty alphabet).
\begin{enumerate}
\item If $\lam^i=\lam^{i-1}$, then $(P_i, T_i)=(P_{i-1}, T_{i-1})$.
\item If $\lam^i \supset  \lam^{i-1}$, then $i=2k$ for some integer
$k \in [n]$. In this case let $P_i=P_{i-1}$ and $T_i$ is
obtained from $T_{i-1}$ by adding the entry $k$ in the square
$\lam^i \setminus \lam^{i-1}$.
\item If $\lam^i \subset \lam^{i-1}$, then $i=2k-1$ for some integer
$k \in [n]$. In this case
let $T_i$ be the unique SYT (on a suitable
alphabet) of shape $\lam^i$ such that $T_{i-1}$ is obtained from $T_i$
by row-inserting some number $j$. Note that $j$ must be less than $k$.
Let $P_i$ be obtained from $P_{i-1}$ by adding the ordered pair $(j, k)$.
\end{enumerate}
It is clear from the above construction that (i)
$P_0 \subseteq P_1 \subseteq \cdots \subseteq P_{2n}$, (ii) for each integer
$i$, it appears at most once as the first component of an ordered pair 
in $P_{2n}$, and appears at most once as the second component of an ordered
pair in $P_{2n}$.   
Let $\psi(V)=(P, T_{2n})$, where
$P$ is the partition on $[n]$ whose standard representation is $P_{2n}$.

Note that 
if an integer $i$ appears in $T_{2n}$, then $P_{2n}$
can not contain any ordered pair $(i,j)$ with $i<j$. It follows that
$i$ is the maximal
element in the block containing it. Hence the content of $T_{2n}$
is a subset of $\max(P)$.

\begin{example}
As an example of the  map $\psi$, let the vacillating tableau be
\[
\eset, \eset, 1, 1, 2, 2, 2, 2, 21, 21, 211, 21, 21, 11, 21.
\]
Then the pairs $(B_i, T_i)$ (where $B_i$ is the pair added to $P_{i-1}$
to obtain $P_i$) are given by
\begin{eqnarray*}
\begin{array}{l|lllllllllllllll}
i     & 0 & 1 & 2 & 3 & 4 & 5 & 6 & 7& 8 & 9  & 10 & 11 & 12 & 13 & 14 \\
\hline
T_i
   &\eset&\eset & 1 & 1 & 12 & 12 & 12 & 12 & 12 & 12 & 12 & 14 & 14 &
   1 & 17\\[-.05in]
   &     &      &   &   &    &    &    &    & 4  & 4  & 4  & 5  & 5  &
   5 & 5 \\[-.05in]
   &     &      &   &   &    &    &    &    &    &    & 5  &    &    &   &   \\
B_i&      &      &   &   &    &    &    &    &    &    &    & (2,6)&  & (4,7)&
\end{array}
\end{eqnarray*}
Hence
\begin{eqnarray*}
T=\begin{array}{l}  1\ 7 \\ [-.05in] 5  \end{array},\quad \quad
P=\mbox{1-26-3-47-5}.
\end{eqnarray*}
\end{example}

The map $\psi$ is bijective since the above construction can be
reversed. Given a pair $(P, T)$, where $P$ is a partition of $[n]$,
and $T$ is an SYT whose content consists of maximal elements of some
blocks of $P$, let $E(P)$ be
the standard representation of $P$, and $T_{2n}=T$. We work our way
backwards from $T_{2n}$, reconstructing the preceding tableaux and
hence the sequence of shapes. If we have the SYT $T_{2k}$ for some
$k \leq n$, we can get the tableaux $T_{2k-1}, T_{2k-2}$ by the following
rules.
\begin{enumerate}
\item
$T_{2k-1}=T_{2k}$ if the integer $k$ does not appear in $T_{2k}$. Otherwise
$T_{2k-1}$ is obtained from $T_{2k}$ by deleting the square containing $k$.
\item
$T_{2k-2}=T_{2k-1}$ if $E(P)$ does not have an edge of the form $(i, k)$.
Otherwise there is a unique $i < k$ such that $(i,k) \in E(P)$.
In that case let  $T_{2k-2}$ be obtained from $T_{2k-1}$ by row-inserting $i$,
or equivalently, $T_{2k-2}=(T_{2k-1} \longleftarrow i)$.
\end{enumerate}

{\bf Proof of Theorem \ref{VT-number}}. \
Part (i)  follows from the bijection $\psi$, where a block of $P$ is
distinguished if its maximal element belongs to $\content(T)$. For
part (ii), simply note that to get a structure counted by $B(n, k)$,
we can partition $[n]$ into two subsets, $S$ and $T$,  and then partition
$S$ into $k$ blocks and put a mark on each block, and partition $T$
arbitrarily. The generating function of $B(n, k)$ then follows from the
well-known generating functions for $S(n,k)$, the Stirling number of the
second kind, and for the Bell number $B(n)$,
\begin{eqnarray*}
\sum_{n \geq k} S(n,k) \frac{x^n}{n!} =\frac{1}{k!} (e^x-1)^k, \qquad
\sum_{n \geq 0} B(n) \frac{x^n}{n!} =\exp(e^x-1). \qquad  \Box
\end{eqnarray*}

\textsc{Remark.} (1). Restricting  to vacillating tableaux of
empty shape, the map $\psi$ provides a bijection between the  set
$\mathcal{V}_\eset^{2n}$ of vacillating tableaux of empty shape
and length $2n$ and the set of partitions of $[n]$. In particular,
$g_\eset(n)$, the cardinality of $\mathcal{V}_\eset^{2n}$, is
equal to the $n$th Bell number $B(n)$.

(2). Note that there is a symmetry between the four types of
 movements in the definition of vacillating tableaux. Thus any walk
 from $\eset$ to $\eset$ in $m+n$ steps can be viewed as a walk from
 $\eset$ to some shape $\lam$ in $n$ steps, then followed by the
 reverse of a walk from $\eset$ to $\lam$ in $m$ steps.  It follows
 that
\begin{eqnarray}\label{bell_number}
\sum_{\lam} g_{\lam}(n) g_\lam(m) = g_\eset(m+n)=B(m+n).
\end{eqnarray}
For the case $m=n=k$, the identity \eqref{bell_number} is proved by 
Halverson and Lowandowski \cite{h-l}, who gave a bijective proof using 
similar procedures as those in $\psi$. 

(3).  The \emph{partition algebra} $\mathfrak{P}_n$ is a
certain semisimple algebra, say over $\mathbb{C}$, whose dimension is
the Bell number $B(n)$ (the number of partitions of $[n]$). (The
algebra $\mathfrak{P}_n$ depends on a parameter $x$ which is
irrelevant here.)  See \cite{h-r,h-l} for a survey of this topic.
Vacillating tableaux are related to irreducible representations of
$\mathfrak{P}_n$ in the same way that SYT of content $[n]$ are related
to irreducible representations of the symmetric group $\sn$.  In
particular, the irreducible representations $I_\lambda$ of
$\mathfrak{P}_n$ are indexed by partitions $\lambda$ for which there
exists a vacillating tableau of shape $\lambda$ and length $2n$, and
$\dim I_n$ is the number of such vacillating tableaux.  This result is
equivalent to \cite[Thm.~2.24(b)]{h-r}, but that paper does not
explicitly define the notion of vacillating tableau.
Combinatorial identities arising from partition algebra and its
subalgebras are discussed in \cite{h-l}, where the authors used the
notion of vacillating tableau after the distribution of   a preliminary 
version of this paper. 

%
%
\section{Crossings and Nestings of Partitions} \label{sec3}
In this section we restrict the map $\psi$ to vacillating tableaux of
empty shape, for which  $\psi$ provides a bijection between
the  set of vacillating tableaux of empty shape and
length $2n$ and the set of partitions of $[n]$.
To make the bijection clear, we restate the inverse map from the set of
partitions to vacillating tableaux.

{\bf The Map $\phi$ from Partitions to Vacillating Tableaux.} \label{phi}
Given a partition $P \in \Pi_n$ with the standard representation,
we construct the sequence of
SYT's, hence the vacillating tableau $\phi(P)$ as follows:
Start from the empty SYT by letting $T_{2n}=\eset$,
read the number $j \in [n]$  one by one from $n$ to 1, and
define $T_{2j-1}$, $T_{2j-2}$ for each $j$.
There are four cases. \\
1.  If  $j$  is the righthand endpoint of an arc $(i,j)$,
but not a lefthand endpoint,  first do nothing, then
insert $i$ (by the RSK algorithm) into the tableau. \\
2. If $j$  is the lefthand endpoint of an arc $(j, k)$,
but not a righthand endpoint, first remove $j$, then do nothing.  \\
3. If $j$ is an isolated point, do nothing twice. \\
4. If $j$  is the righthand endpoint of an arc $(i,j)$,
and the lefthand endpoint of another arc $(j,k)$, then
delete $j$ first, and then insert $i$. \\
The vacillating tableau $\phi(P)$
  is the sequences of shapes of the above  SYT's.

\begin{example} \label{VT1_SYT}
Let $P$ be the partition 1457-26-3 of $[7]$.
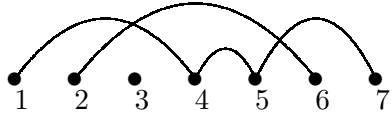
\begin{figure}[ht]
\begin{center}
\begin{picture}(120,40)
\setlength{\unitlength}{4mm} \multiput(0,0)(2,0){7}{\circle*{0.4}}
\qbezier(0,0)(3,4)(6,0) \qbezier(2,0)(6,5)(10,0)
\qbezier(6,0)(7,2)(8,0) \qbezier(8,0)(10,4)(12,0)
\put(0,-1){$1$}
\put(2,-1){$2$}\put(4,-1){$3$}\put(6,-1){$4$}\put(8,-1){$5$}\put(10,-1){$6$}\put(12,-1){$7$}
\end{picture}
\end{center}
\caption{The standard representation of the  partition 1457-26-3.}
\label{Ex_1}
\end{figure}

Starting from $\eset$ on the right, go from 7 to 1, the seven steps are
(1) do nothing, then insert 5,
(2) do nothing, then insert 2,
(3) delete 5 and insert 4,
(4) delete 4 and insert 1,
(5) do nothing twice,
(6) remove 2 then do nothing,
(7) remove 1 then do nothing.
 Hence the corresponding SYT's, constructed from right to left, are
\begin{eqnarray*}
\begin{array}{ccccccccccccccc}
\eset & \eset & 1 & 1 & 1 & 1 & 1 & 2 & 24 & 2 & 2 & 5 & 5 & \eset
      &\eset \\[-.05in]
      &       &   &   & 2 & 2 & 2 &   &    &   & 5 &   &   &       &
\end{array}
\end{eqnarray*}
The vacillating tableau is
\begin{eqnarray*}
\eset, \eset, 1, 1,11, 11, 11, 1, 2, 1, 11, 1, 1, \eset, \eset.
\end{eqnarray*}
\end{example}

The relation between $\crr(P), \nee(P)$ and the vacillating tableau is
given in the next theorem.

\begin{theorem} \label{CN}
Let $P \in \Pi_n$ and $\phi(P)=(\eset=\lam^0, \lam^1, \dots, \lam^{2n}=\eset)$.
Then $\crr(P)$ is the most number of rows in any $\lam^i$, and $\nee(P)$ is
the most number of columns in any $\lam^i$.
\end{theorem}
\noindent {\bf Proof}.
 We prove Theorem~\ref{CN} in four steps.
First, we interpret a $k$-crossing/$k$-nesting of $P$ in terms of
entries of SYT's  $T_i$ in $\phi(P)$. Then, we associate to each SYT $T_i$
a sequence $\sigma_i$ whose terms are entries of $T_i$.  We prove that
$T_i$ is the insertion tableau of $\sigma_i$ under the RSK algorithm,
and  apply Schensted's theorem to conclude the proof.

{\bf Step 1}.
Let $T(P)=(T_0, T_1, \dots, T_{2n})$ be the sequence
of  SYT's associated to the vacillating tableau $\phi(P)$.
By the construction of $\psi$ and $\phi$, a pair $(i, j)$ is an arc in
the standard representation of $P$ if and only if $i$ is an entry in the
SYT's $T_{2i}, T_{2i+1}, \dots, T_{2j-2}$. We say that the the integer $i$
is \emph{added} to $T(P)$ at step $i$ and \emph{leaves}  at step $j$.

First we prove that the arcs  $(i_1, j_1), \dots,(i_k, j_k)$ form
a $k$-crossing of $P$ if and only if there exists a tableau $T_i$
in $T(P)$ such that the integers $i_1, i_2, \dots, i_k \in
\content(T_i)$, and $i_1, i_2, \dots, i_k$ leave  $T(P)$ in
increasing order according to their numerical values. Given a
$k$-crossing $((i_1, j_1), \dots,(i_k, j_k))$ of $P$, where $i_r <
j_r$ for $1 \leq r \leq k$ and $i_1< i_2< \cdots < i_k < j_1 < j_2
< \cdots < j_k$, the integer $i_r$ is added to $T(P)$ at step
$i_r$ and leaves at step $j_r$. Hence all $i_r$ are in
$T_{2j_1-2}$, and they leave $T(P)$ in increasing order. The
converse is also true: if there are $k$ integers $i_1 < i_2 <
\dots < i_k$ all appearing in the same tableau at some steps, and
then leave in increasing order, say at steps $j_1 < j_2 < \dots <
j_k$, then $i_k < j_1$ and the pairs $((i_1, j_1), \dots,(i_k,
j_k)) \in P$ form a $k$-crossing. By a similar argument  arcs
 $(i_1, j_1), \dots,(i_k, j_k)$ form a  $k$-nesting of $P$
if and only if there exists a tableaux $T_i$ in $T(P)$ such that
the integers $i_1, i_2, \dots, i_k \in \content(T_i)$, and $i_1, i_2, \dots,
i_k$ leave $T(P)$ in decreasing order.

{\bf Step 2}.
For each $T_i \in T(P)$, we define a permutation $\sigma_i$ of
$\content(T_i)$  (backward) recursively as follows.
Let $\sigma_{2n}$ be the empty sequence.
(1) If $T_i=T_{i-1}$, then $\sigma_{i-1}=\sigma_i$.
(2) If $T_{i-1}$ is obtained from $T_i$ by row-inserting some number $j$,
then $\sigma_{i-1}=\sigma_ij$, the juxtaposition of $\sigma_i$ and $j$.
(3) If $T_i$ is obtained from $T_{i-1}$ by adding the entry $i/2$,
(where $i$ must be even),
then $\sigma_{i-1}$ is obtained from $\sigma_i$ by deleting the number $i/2$.
Note that in the last case, $i/2$ must be the largest entry in $\sigma_i$.

Clearly $\sigma_i$ is a permutation of the entries in $\content(T_i)$.
If $\sigma_i=w_1w_2\cdots w_r$, then the entries of $\content(T_i)$
leave $T(P)$ in the order $w_r, \dots, w_2, w_1$.

{\bf Step 3.}
{\bf Claim}:
If $\sigma_i \stackrel{\text{RSK}}{\longmapsto} (A_i, B_i)$, then $A_i=T_i$.

We prove the claim by backward induction. The case $i=2n$ is trivial
as both $A_{2n}$ and $T_{2n}$ are the empty SYT.  Assume the claim is
true for some $i$, $ 1 \leq i \leq 2n$.
We prove that the claim holds for $i-1$.

If $T_{i-1}=T_i$,  the claim holds by the inductive hypothesis.
If  $T_{i-1}$ is obtained from
$T_i$ by inserting the number $j$, then the claim holds  by
the definition of the RSK algorithm. It is only left to consider the
case that $T_{i-1}$ is obtained from $T_{i}$ by removing the entry
$j=i/2$.

Let us write $\sigma_i$ as
$u_1u_2\cdots u_s j v_1\cdots v_t$, and $\sigma_{i-1}$ as
$u_1u_2\cdots u_s v_1 \cdots v_t$, where $j > u_1, \dots,u_s, v_1,
\dots, v_t$.  We need to show that the insertion tableau of
$\sigma_{i-1}$ is the same as the insertion tableau of $\sigma_i$
deleting the entry $j$, i.e., $A_{i-1}=A_i\!\setminus\!
\{j\}$.  Proof by induction on $t$.  If $t=0$ then it is true by the
RSK algorithm that $A_i$ is obtained from $A_{i-1}$ by
adding $j$ at the end of the first row. Assume it is true for $t-1$,
i.e., $A(u_1\cdots u_s v_1\cdots v_{t-1})=A(u_1\cdots u_s j v_1\cdots
v_{t-1}) \setminus \{j\}$.  Note that in $A(u_1\cdots u_s j v_1\cdots
v_{t-1})$, if $j$ is in position $(x, y)$, then there is no element in
positions $(x, y+1)$ or $(x+1, y)$.  Now we insert the entry $v_t$ by
the RSK algorithm.  Consider the insertion path $I=I(A(u_1\cdots
u_s j v_1\cdots v_{t-1})\longleftarrow v_t ). $ If $j$ does not appear
on this path, then we would have the exact same insertion path when
inserting $v_t$ into $A(u_1\cdots u_sv_1\cdots v_{t-1})$. This
insertion path results in the same change to $A(u_1\cdots u_s j
v_1\cdots v_{t-1})$ and $A(u_1\cdots u_s v_1\cdots v_{t-1})$, which
does not touch the position $(x,y)$ of $j$. So $A(u_1\cdots u_s
v_1\cdots v_{t-1}v_t )= A(u_1\cdots u_s j v_1\cdots
v_{t-1}v_t)\setminus \{j\}$.  On the other hand, if $j$ appears in the
insertion path $I$, i.e., $(x,y) \in I$, then since $j$ is the largest
element, it must be bumped into the $(x+1)$-th row, and become the
last entry in the $(x+1)$-th row without bumping any number further.
Then the insertion path of $v_t$ into $A(u_1\cdots u_sv_1\cdots
v_{t-1})$ is $I$ minus the last position $\{(x+1, *)\}$, and again we
have $A(u_1\cdots u_s v_1\cdots v_{t-1}v_t )= A(u_1\cdots u_s j
v_1\cdots v_{t-1}v_t)\setminus \{j\}$.  This finishes the proof of the
claim.

{\bf Step 4}. We shall need the following theorem of Schensted
\cite{Sch61}\cite[Thms.~7.23.13, 7.23.17]{Stanley99}, which gives the
basic connection between the RSK algorithm and the increasing and
decreasing subsequences.

\noindent

\begin{quote}
\noindent {\bf Schensted's Theorem}
Let $\sigma$ be a sequence of integers whose terms are distinct. Assume
$\sigma \stackrel{RSK}{\longmapsto} (A, B)$, where $A$ and $B$ are SYT's of
the shape  $\lam$.
Then the length of the longest increasing subsequence of $\sigma$ is
$\lam_1$ (the number of columns of $\lam$), and the length of the longest
decreasing subsequence is $\lam'_1$ (the number of rows of $\lam$).
\end{quote}

Now we are ready to prove Theorem~\ref{CN}. By Steps 1 and 2,  a
partition $P$ has a $k$-crossing if and only if there exists $i$
such that $\sigma_i$ has a decreasing subsequence of length $k$.
The claim in Step 3 implies that the shape of the sequence
$\sigma_i$ is exactly the diagram of the $i$-th partition $\lam^i$
in the vacillating tableau $\phi(P)$. By Schensted's Theorem,
$\sigma_i$ has a decreasing subsequence of length $k$ if and only
if the partition $\lam^i$ in $\phi(P)$ has at least $k$ rows. This
proves the statement for $\crr(P)$ in  Theorem~\ref{CN}. The
statement  for $\nee(P)$ is proved similarly. ~~~$\Box$

\vspace{.3cm} The symmetric joint distribution of statistics
$\crr(P)$ and $\nee(P)$ over $P_n(S, T)$ follows immediately from
Theorem~\ref{CN}.

\noindent {\bf Proof of Theorem 1.}  \\

From Theorem~\ref{CN}, a partition $P \in \Pi_n$ has $\crr(P)=k$ and
$\nee(P)=j$ if and only if for the partitions $\{ \lam^i\}_{i=0}^{2n}$
of the vacillating tableau $\phi(P)$, the maximal number of rows of
the diagram of any $\lam^i$ is $k$, and the maximal number of columns
of the diagram of any $\lam^i$ is $j$.  Let $\tau$ be the involution
defined on the set ${\cal V}_\emptyset^{2n}$ by taking the conjugate
to each partition $\lam^i$.  For $i \in [n]$, $i \in \min(P)$
(resp. $\max(P)$) if and only if $\lambda^{2i-1}=\lambda^{2i-2}$ 
and $\lambda^{2i} \setminus \lambda^{2i-1} = \Box$, 
(reps. $\lambda^{2i-2} \setminus  \lambda^{2i-1}=\Box$ and
$\lambda^{2i}=\lambda^{2i-1}$).  
Since $\tau$ preserves $\min(P)$ and
$\max(P)$, it induces an involution on $P_n(S, T)$ which exchanges
the statistics $\crr(P)$ and $\nee(P)$. This proves Theorem
\ref{thm1}.
~~~$\Box$ \vspace{.3cm}

Let $\lam=(\lam_1, \lam_2, \dots )$ be the shape of a sequence $w$ of
distinct integers.
Schensted's Theorem provides a combinatorial interpretation of the terms
$\lam_1$ and $\lam'_1$: they are the length of the longest increasing
and decreasing subsequences of  $w$.
In \cite{Greene74} C. Greene extended Schensted's Theorem by giving an
interpretation of  the
rest of the diagram of $\lam=(\lam_1, \lam_2, \dots )$.

Assume $w$ is a sequence of length $n$. For each $k \leq n$, let
$d_k(w)$ denote the length of the longest subsequence of $w$ which
has no increasing subsequences of length $k+1$. It can be shown easily that any
such sequence is obtained by taking the union of $k$ decreasing subsequences.
Similarly, define $a_k(w)$ to be the length of the longest subsequence
consisting of $k$ ascending subsequences.

\begin{theorem}[Greene]
For each $k \leq n$,
\begin{eqnarray*}
a_k(w)&=&\lam_1+\lam_2+\cdots +\lam_k, \\
d_k(w)& =& \lam'_1 +\lam'_2 +\cdots + \lam'_k,
\end{eqnarray*}
where $\lam'=(\lam_1', \lam_2', \dots)$ is the conjugate of $\lam$.
~~~~$\Box$
\end{theorem}

We may consider the analogue of Greene's Theorem for set partitions.
Let  $P \in \Pi_n$  with the standard
representation $\{(i_1, j_1), (i_2, j_2), \dots, (i_t, j_t)\}$,
where $i_r < j_r$ for $1 \leq r \leq t$.
 Let $e_r=(i_r, j_r)$. We define
 the \emph{crossing graph} $\mathrm{Cr}(P)$ of $P$ as follows. The vertex set
of $\mathrm{Cr}(P)$ is $\{e_1, e_2, \dots, e_t\}$.  Two arcs $e_r$ and
$e_s$ are adjacent if and only if the edges $e_r$ and $e_s$ are
crossing, that is, $i_r < i_s < j_r < j_s$.
Clearly a $k$-crossing of $P$ corresponds to a
$k$-clique of $\mathrm{Cr}(P)$. Let $\crr_r(P)$ be the maximal number
of vertices in a union of $r$ cliques of $\mathrm{Cr}(P)$. In other
words, $\crr_r(P)$ is the maximal number of arcs in a union of $r$
crossings of $P$.  Similarly, let $\mathrm{Ne}(P)$ be the graph
defined on the vertex set $\{e_1, \dots, e_t\}$ where two arcs $e_r$
and $e_s$ are adjacent if and only if $i_r < i_s < j_s < j_r$. Let
$\nee_r(P)$ be the maximal number of vertices in a union of $r$
cliques of $\mathrm{Ne}(P)$. In other words, $\nee_r(P)$ is the maximal number
of arcs in a union of $r$ nestings of $P$.

\begin{prop} \label{Greene1}
Let $P=((i_1, j_1), (i_2, j_2), \dots, (i_t, j_t))$ be the
standard representation of a partition of $[n]$, where $i_r < j_r$
for all $1 \leq r \leq t$ and $j_1 < j_2 < \cdots< j_t$. Let
$\alpha(P)$ be the sequence $i_1i_2\cdots i_t$. Then  there is a
one-to-one correspondence between the set of nestings of $P$ and
the set of decreasing subsequences of $\alpha(P)$.
\end{prop}
\noindent {\bf Proof}.
Let $\phi(P)$ be the vacillating tableau corresponding to $P$,
and $T(P)$ the sequence of SYT's constructed  in the
bijection.
Then $\alpha(P)$ records  the order in which the entries of
 $T_i$'s leave $T(P)$.
Let $\sigma = i_t \cdots i_2 i_1$ be the reverse of $\alpha(P)$,
and $\{\sigma_i : 1 \leq i \leq 2n\}$ the permutation of $\content(T_i)$
defined in Step 2 of the proof of Theorem \ref{CN}.
Then the $\sigma_i$'s  are subsequences of $\sigma$.

From Steps 1 and 2 of the proof of Theorem~\ref{CN}, nestings
of $P$ are represented by the increasing subsequences of $\sigma_i$,
$1 \leq i \leq 2n$, and hence by the increasing subsequences of
$\sigma$. Conversely, let $i_{r_1} <i_{r_2}
 < \cdots < i_{r_t}$ be an increasing subsequence of $\sigma$.
Being a subsequence of $\sigma$ means that its terms leave $T(P)$ in
reverse order,
so $i_{r_t}$ leaves first in step $j_{r_t}$. Thus
all the entries $i_{r_1}, \dots, i_{r_t}$ appear in the SYT's $T_i$ with
$2i_{r_t} \leq i \leq 2j_{r_t}-2$. Therefore
$i_{r_1}i_{r_2}\cdots i_{r_t}$  is  also an increasing subsequence
of  $\sigma_i$, for $2i_{r_t} \leq i \leq 2j_{r_t}-2$.
~~~ $\Box$

Combining Proposition~\ref{Greene1} and Greene's Theorem, we have
the following corollary describing $\nee_t(P)$.
\begin{cor} \label{cor9}
Let $P$ and $\alpha(P)$ be as in Proposition \ref{Greene1}.
 Then
\[
\nee_r(P)=\lam'_1+\lam'_2+\cdots
+\lam'_r,
\]
where $\lam$ is the shape of $\alpha(P)$,
and $\lam'$ is the conjugate of $\lam$. ~~~~$\Box$
\end{cor}

The situation for $\crr_r(P)$ is more complicated. We don't have a
result similar to Proposition~\ref{Greene1}. Any crossing of $P$
uniquely corresponds to an increasing subsequence of $\alpha(P)$. But
the converse is not true. An increasing subsequence $i_{r_1}i_{r_2}
\cdots i_{r_t}$ corresponds to a $t$-crossing of $P$ only if we have
the additional condition $i_{r_t} < j_{r_1}$. It would be interesting
to get a result for $\crr_r(P)$ analogous to Corollary~\ref{cor9}.

To conclude this section we discuss the enumeration of noncrossing
partitions.
The following theorem is a direct corollary of
the bijection between vacillating tableaux and partitions.
\begin{theorem} \label{lattice_P}
Let $\epsilon_i$ denote the $i$th unit coordinate vector in
$\rr^{k-1}$. The number of $k$-noncrossing partitions of $[n]$
 equals the number of closed lattice walks in the region
\[
V_k=\{(a_1, a_2, \dots, a_{k-1}): a_1 \geq a_2 \geq \cdots \geq
a_{k-1} \geq 0, a_i \in \zz \}
\]
from the origin to itself of length $2n$ with steps
$\pm\epsilon_i$ or $(0,0,\dots,0)$, with the property that the
walk goes backwards (i.e., with step $-\epsilon_i$) or stands
still (i.e., with step $(0,0,\dots,0)$) after an even number of
steps, and  goes forwards (i.e., with step $+\epsilon_i$) or
stands still after an odd number of steps.
~~~~$\Box$
\end{theorem}

Recall that a partition $P \in \Pi_n$ is \emph{$k$-noncrossing} if
$\crr(P) < k$, and is \emph{$k$-nonnesting} if $\nee(P) <k$.
A partition $P$ has
 no  $k$-crossings and no $j$-nestings
if and only if for all the partitions $\lambda^i$ of $\phi(P)$,
the diagram
fits into a $(k-1)\times (j-1)$ rectangle. Taking the conjugate of each
partition, we get  bijective proofs of Corollaries \ref{ncn} and
\ref{nc_nn}.

Theorem \ref{thm1} asserts the symmetric distribution of $\crr(P)$ and
$\nee(P)$ over $P_n(S, T)$, for all $S, T \subseteq [n]$ with $|S|=|T|$.
Not every $P_n(S,T)$ is nonempty. A set
$P_n(S, T)$ is nonempty if and only if for all $i \in [n]$,
$|S \cap [i]| \geq |T \cap [i]|$.  Another way to describe the nonempty
$P_n(S, T)$
is to use lattice paths. Associate to each pair $(S,T)$ a lattice path
$L(S, T)$  with steps $(1, 1)$, $(1, -1)$ and $(1, 0)$:
start from $(0,0)$, read the integers $i$ from $1$ to $n$ one by one, and
move  two steps for each $i$. \\
\mbox{} \hspace{2em} 1. If $i \in S \cap T$, move $(1, 0)$ twice.\\
\mbox{} \hspace{2em} 2. If $i \in S\setminus T$, move $(1,0)$ then $(1,1)$.\\
\mbox{} \hspace{2em} 3. If $i \in T\setminus S$, move $(1,-1)$ then $(1,0)$.\\
\mbox{} \hspace{2em} 4. If $i \notin S\cup T$, move $(1,-1)$ then $(1,1)$. \\
This defines a lattice path $L(S, T)$ from $(0,0)$ to $(2n, 0)$,
Conversely, the path uniquely determines $(S, T)$.
Then $P_n(S, T) $ is nonempty if and only if the lattice path $L(S,
T)$ is a Motzkin path, i.e., never goes below the $x$-axis.

There are existing notions of noncrossing partitions and nonnesting
partitions, e.g., \cite[Ex.6.19]{Stanley99}.
A  \emph{noncrossing partition of $[n]$} is a partition of $[n]$
in which no two blocks ``cross'' each other, i.e.,
if $a < b < c <d$ and $a, c$ belong to a  block $B$ and $b, d$ to
another block  $B'$,
then $B=B'$. A \emph{nonnesting partition of $[n]$} is a partition
of $[n]$ such that if $a, e$ appear in a block $B$ and $b, d$ appear
in a different block $B'$  where $a < b < d< e$, then there is a $c
\in B$ satisfying $b < c < d$.

It is easy to see that $P$ is a noncrossing partition
if and only if the standard representation of $P$ has no 2-crossing, and
$P$ is a nonnesting partition if and only if the standard representation of $P$
has no 2-nesting.   Hence
the vacillating tableau correspondence, in the case of 1-row/column
tableaux, gives  bijections
between noncrossing partitions of $[n]$, nonnesting partitions
of $[n]$, (both are counted by Catalan numbers) and
sequences $0=a_0, a_1, ..., a_{2n}=0$ of nonnegative integers such that
$a_{2i+1} = a_{2i}$ or $a_{2i}-1$, and $a_{2i} = a_{2i-1}$ or
$a_{2i-1}+1$. These sequences $a_0, ..., a_{2n}$ give a new combinatorial
interpretation of Catalan numbers.

Replacing a term $a_{i+1} = a_i + 1$ with a step $(1,1)$, a
term $a_{i+1} = a_i - 1$ with a step $(1,-1)$, and a term $a_{i+1} = a_i$
with a step $(1,0)$, we get a Motzkin path, so we also have a bijection
between noncrossing/nonnesting partitions and certain Motzkin paths.
The Motzkin paths are exactly the ones defined as $L(S, T)$, where
$S=\min(P)$ and $T=\max(P)$.
Conversely, given a Motzkin path of the form $L(S, T)$,
we can recover uniquely a noncrossing partition and a nonnesting
partition. Write the path as  $\{(i, a_i): 0\leq i \leq 2n\}$.
Let $A=[n] \setminus T$ and $B=[n]\setminus S$.
Clearly $|A|=|B|$.
Assume
$A=\{ i_1, i_2, \dots, i_t\}_<$, and $B=\{ j_1, j_2, \dots, j_t\}_<$
where elements are listed in increasing order. Then to get
the standard representation of the noncrossing partition, pair each
$j_r$ with $\max\{i_s \in A: i_s < j_r, a_{2i_s}=a_{2j_r-2}\}$.
To get the standard representation of the nonnesting partition,
pair each $j_r $ with $i_r$, for $1\leq r \leq t$.

{\sc Remark.}
In our definition, a  $k$-crossing is defined as a set of
$k$ mutually crossing arcs in the standard representation of the
partition.
There exist some other definitions. For example, in \cite{Klazar03}
M. Klazar defined \emph{3-noncrossing partition} as
a partition $P$ which does not have 3 mutually crossing blocks.
It can be seen that $P$ is 3-noncrossing in Klazar's sense
if and only if there do not
exist 6 elements
$a_1 < b_1 < c_1 < a_2 < b_2 < c_2$ in $[n]$  such that $a_1, a_2 \in A$,
$b_1, b_2 \in B$,  $c_1, c_2 \in C$, and
$A, B, C$ are three distinct blocks of $P$.

Klazar's definition of 3-noncrossing partitions is different from ours.
For example, let $P$ be the partition 15-246-37 of $[7]$, with
standard representation as follows:

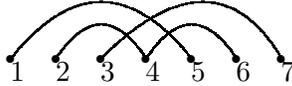
\begin{center}
\begin{picture}(120,30)
\setlength{\unitlength}{3mm}
\put(0,0){\circle*{0.4}}
\put(2,0){\circle*{0.4}}
\put(4,0){\circle*{0.4}}
\put(6,0){\circle*{0.4}}
\put(8,0){\circle*{0.4}}
\put(10,0){\circle*{0.4}}
\put(12,0){\circle*{0.4}}
\qbezier(0,0)(4,5)(8,0)
\qbezier(2,0)(4,3)(6,0)
\qbezier(4,0)(9,5)(12,0)
\qbezier(6,0)(8,3)(10,0)
\put(0,-1){$1$}
\put(2,-1){$2$}
\put(4,-1){$3$}
\put(6,-1){$4$}
\put(8,-1){$5$}
\put(10,-1){$6$}
\put(12,-1){$7$}
\end{picture}
\end{center}

\medskip
According to Klazar's definition, $P$
has a 3-crossing, since we have $1<2<3<5<6<7$ and $\{1,5\}$, $\{2,6\}$
and $\{3,7\}$ belong to three different blocks, respectively.
On the other hand, $P$ has no 3-crossing on our sense.

For general $k$, these three notions of $k$-noncrossing partitions, i.e.,
(1) no $k$-crossing in the standard representation of $P$,
(2) no $k$ mutually crossing arcs in  distinct blocks of $P$,
and (3) no $k$ mutually crossing blocks, are all different, with
the first being the weakest, and the third  the strongest.

%
%
\section{ A Variant: Partitions and Hesitating Tableaux}
We may also consider the \emph{enhanced}  crossing/nesting of a partition,
by taking isolated points into consideration. For a partition $P$ of $[n]$,
 let the \emph{enhanced representation} of  $P$ be
the union of  the standard representation of $P$  and the
 loops $\{(i,i): i \text{ is an isolated point of $P$}\}$.
An \emph{enhanced $k$-crossing} of $P$
is a set of $k$ edges  $(i_1, j_1), (i_2, j_2), \dots, (i_k, j_k)$ of
the enhanced  representation of $P$
such that $i_1 < i_2 < \cdots < i_k \leq j_1 < j_2 < \cdots < j_k$.
In particular, two arcs of the form $(i,j)$ and $(j,l)$
with $i <j <l$ are viewed as  crossing.
Similarly, an \emph{enhanced $k$-nesting} of $P$  is a set of
$k$ edges $(i_1, j_1), (i_2, j_2), \dots, (i_k, j_k)$ of the
enhanced representation of $P$  such that
$i_1 < i_2 < \cdots < i_k \leq j_k < \cdots < j_2 < j_1$.
In particular, an edge $(i,k)$ and an isolated point $j$ with $i < j<k$
form an enhanced 2-nesting.

Let $\overline{\crr}(P)$ be the size of the largest enhanced crossing,
and $\overline{\nee}(P)$ the size of the largest enhanced nesting.
Using a variant form of vacillating tableau, we obtain
again a symmetric joint distribution of
the statistics $\overline{\crr}(P)$ and $\overline{\nee}(P)$.

The variant tableau is a
\emph{hesitating tableau} of shape $\emptyset$ and length $2n$,
which is a path on the Hasse diagram of Young's lattice from $\eset$
to $\eset$ where each step consists of a pair of moves, where the pair
is either (i)  doing nothing then
adding a square, (ii) removing a square then doing nothing, or (iii)
adding a square and then removing a square.

\begin{example} \label{GOT}
There are 5 hesitating  tableaux of shape $\emptyset$ and length
$6$. They are
\begin{eqnarray} \label{GOT_3}
\begin{array}{ccccccc}
\emptyset  & 1 & \eset  &   1     & \eset & 1    & \eset \\
\eset      & 1 & \eset  &   \eset & 1     & \eset& \eset \\
\eset      &  \eset & 1 &  11     &  1    & \eset&  \eset  \\
\eset      &  \eset & 1 &  \eset  & \eset  & 1   & \eset \\
\eset      &  \eset & 1 & 2      & 1 & \eset &\eset
\end{array}
\end{eqnarray}
\end{example}

To see the equivalence with vacillating tableaux, let
$U$ be the operator that takes a shape to the sum of all shapes that
cover it in Young's lattice (i.e., by adding a square), and similarly
$D$ takes a shape to the sum of all shapes that it covers in Young's
lattice (i.e., by deleting a square). Then, as is well-known,
$DU-UD=I$ (the identity operator). See,
e.g.\ \cite{Stanley88}\cite[Exer.~7.24]{Stanley99}. It follows that
\begin{eqnarray} \label{opr}
(U+I)(D+I)=DU+ID+UI.
\end{eqnarray}
Iterating the left-hand side generates vacillating tableaux, and
iterating the right-hand side gives the hesitating tableaux defined
above.

A bijective map between partitions of $[n]$ and hesitating
tableaux of empty shape has been given by Korn \cite{Korn},
based on growth diagrams. Here by modifying the map $\phi$ defined
in Section~\ref{sec3} we get a more direct bijection $\bar \phi$
between partitions and hesitating tableaux of empty shape, which
leads to the symmetric joint distribution of $\overline{\crr}(P)$
and $\overline{\nee}(P)$.  The construction and proofs are very
similar to the ones given in Sections~\ref{sec2} and \ref{sec3},
and hence are omitted here. We will only state the definition of
the map $\bar \phi$ from partitions to hesitating tableaux, to be
compared with the map $\phi$ in Section~\ref{sec3}.

{\bf The Bijection $\bar \phi$ from Partitions to Hesitating Tableaux}. \\
Given a partition $P \in\Pi_n$  with the enhanced representation,
we construct the sequence of
SYT's, and hence the hesitating tableau $\bar \phi(P)$, as follows:
Start from the empty SYT by letting $T_{2n}=\eset$,
read the numbers $j \in [n]$  one by one from $n$ to $1$, and
 define two SYT's $T_
{2j-1}$, $T_{2j-2}$ for each $j$.
When $j$ is a lefthand endpoint only, or a righthand endpoint only,
the construction is identical to that of the map $\phi$. Otherwise, \\
1. If $j$ is an isolated point, first insert $j$, then delete $j$.\\
2. If $j$  is the righthand endpoint of an arc $(i,j)$,
and the lefthand endpoint of another arc $(j,k)$, then
insert $i$ first, and then delete $j$.

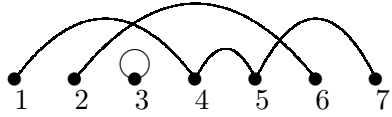
\begin{figure}[ht]
\begin{center}
\begin{picture}(120,40)
\setlength{\unitlength}{4mm} \multiput(0,0)(2,0){7}{\circle*{0.4}}
\qbezier(0,0)(3,4)(6,0) \qbezier(2,0)(6,5)(10,0)
\qbezier(6,0)(7,2)(8,0) \qbezier(8,0)(10,4)(12,0)
\put(4,0.5){\circle{1}} \put(0,-1){$1$}
\put(2,-1){$2$}\put(4,-1){$3$}\put(6,-1){$4$}\put(8,-1){$5$}\put(10,-1){$6$}\put(12,-1){$7$}
\end{picture}
\end{center}
\caption{The enhanced representation of the partition 1457-26-3.}
\label{Ex_2}
\end{figure}

\begin{example} \label{GOT_SYT}
For the partition  1457-26-3 of $[7]$ in Figure \ref{Ex_2},
the corresponding SYT's are
\begin{eqnarray*}
\begin{array}{lllllllllllllll}
\eset & \eset & 1 & 1 & 1 & 13 & 1 & 14 & 24 & 24 & 2 & 5 & 5 & \eset
      & \eset\\[-.05in]
      &       &   &   & 2 & 2  & 2 & 2  &    & 5  & 5 &   &   &       &
\end{array}
\end{eqnarray*}
The hesitating tableau $\bar \phi(P)$ is
\[
\eset,  \eset,  1, 1, 11, 21, 11, 21, 2, 21, 11, 1, 1, \eset, \eset.
\]
\end{example}

The conjugation of shapes does not preserve $\min(P)$ or $\max(P)$.
Instead, it preserves $\min(P)\!\setminus\! \max(P)$, and $\max(P)\!
\setminus\! \min(P)$. Let $S, T$ be disjoint subsets of $[n]$ with the
same cardinality,
$\bar P_n(S, T)=\{ P \in \Pi_n: \min(P)\!\setminus\! \max(P)=S,
\max(P) \!\setminus\! \min(P)=T\}$, and
$\bar f_{n, S, T}(i,j)=\#\{ P \in \bar P_n(S, T): \overline{\crr}(P)=i,
\overline{\nee}(P)=j\}$.
\begin{theorem}
We have
\[
\bar f_{n,S, T}(i,j)=\bar f_{n, S, T} (j, i). \qquad \Box
\]
\end{theorem}
As a consequence, Corollaries \ref{ncn} and \ref{nc_nn} remain valid
if we define $k$-noncrossing (or $k$-nonnesting) by $\bar \crr(P)<k$
(or $\bar \nee(P)<k$).

\textsc{Remark.} As done for vacillating tableaux, one can extend the
 definition of hesitating tableaux by considering moves from $\eset$ to
 $\lam$, and denote by $f_{\lam}(n)$ the number of such hesitating
 tableaux of length $2n$. Identity \eqref{opr} implies that
 $f_{\lam}(n)=g_{\lam}(n)$, the number of vacillating tableaux from
 $\eset $ to $\lam$. It follows that $f_{\eset}(n)=B(n)$, and
\[
\sum_\lam f_\lam(n) f_\lam(m)=B(m+n),
\]
where $B(n)$ is the $n$th Bell number. For further discussion of the
number  $f_\lam(n)$, see  \cite[Problem 33]{Stanley_web} (version of
17 August 2004).

%
%
\section{Enumeration of $k$-Noncrossing  Matchings}
Restricting Theorem~\ref{thm1} to disjoint subsets $(S, T)$ of $[n]$,
where $n=2m$ and $|S|=|T|=m$, we get the symmetric joint distribution of the
crossing number and nesting number for matchings, as stated in
Corollary \ref{matching} in Section 1.

In a complete matching, an integer is either a left endpoint or a
right endpoint in the standard representation.  In applying the map
$\phi$ to complete matchings on $[2m]$, if we remove all steps which
do nothing, we obtain a sequence $\emptyset=\lam^0$, $\lam^1, \dots,
\lam^{2m} =\eset$ of partitions such that for all $1 \leq i \leq 2m$,
the diagram of $\lam^i$ is obtained from that of $\lam^{i-1}$ by
either adding one square or removing one square. Such a sequence is
called an \emph{oscillating tableau} (or \emph{up-down tableau}) of
empty shape and length $2m$.  Thus we get a bijection between complete
matchings on $[2m]$ and oscillating tableaux of empty shape and
length $2m$.  This bijection was originally constructed by the fourth
author, and then extended by Sundaram \cite{Sundaram90} to arbitrary
shapes to give a combinatorial proof of the Cauchy identity for the
symplectic group Sp$(2m)$.  The explicit description of the bijection
has appeared in \cite{Sundaram90} and was included in \cite[Exercise
7.24]{Stanley99}. Oscillating tableaux first appeared (though not with
that name) in \cite{berele}.

Recall that the ordinary RSK algorithm gives a bijection between the 
symmetric group  $\sm$ and pairs $(P, Q)$ of SYTs of the same shape 
$\lambda \vdash m$. This result and the Schensted's theorem  
can be viewed as a special case of what we do. 
Explicitly, identify an SYT $T$ of shape $\lambda$ and content $[m]$  
with a sequence $\emptyset =\lam^0, \lam^1, \dots, \lam^m=\lam$ of integer
partitions where $\lam^i$ is the shape of  the SYT $T^i$ obtained from 
$T$ by deleting all entries $\{j: j >i\}$. Let $w$ be a permutation of
$[m]$, and form the matching $M_w$ on $[2m]$ with arcs between $w(i)$ 
and $2m-i+1$. We get an oscillating tableau $O_w$ that increases to the 
shape $\lam \vdash m$ and then decreases to the empty shape. Assume $w 
\stackrel{\text{RSK}}{\longmapsto}(A(w), B(w))$, where $A(w)$ is the
(row)-insertion tableau and $B(w)$ the recording tableau. Then $A(w)$ 
is given by the first $m$ steps of $O_w$, and $B(w)$ 
the reverse of the last $m$ steps. 
The size of the largest crossing (resp. nesting) of $M_w$ is exactly
the length of the longest decreasing (resp. increasing) subsequence of $w$. 

\begin{example} 
Let $w=231$. Then 
\begin{eqnarray*}
A(w)=\begin{array}{cc}  1 & 3 \\
                        2 & \end{array}, 
\qquad 
B(w)=\begin{array}{cc}  1 & 2 \\
                        3 & \end{array}. 
\end{eqnarray*} 
The matching $M_w$ and the corresponding oscillating tableau are given 
below. 
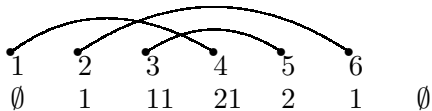
\begin{center}
\begin{picture}(120,30)
\setlength{\unitlength}{3mm}
\put(0,1){\circle*{0.4}}
\put(3,1){\circle*{0.4}}
\put(6,1){\circle*{0.4}}
\put(9,1){\circle*{0.4}}
\put(12,1){\circle*{0.4}}
\put(15,1){\circle*{0.4}}
\qbezier(0,1)(4,4)(9,1)
\qbezier(3,1)(9,5)(15,1)
\qbezier(6,1)(9,3)(12,1)
\put(0,0){$1$}
\put(3,0){$2$}
\put(6,0){$3$}
\put(9,0){$4$}
\put(12,0){$5$}
\put(15,0){$6$}
\put(0,-1.5){$\emptyset$}
\put(3,-1.5){$1$}
\put(6,-1.5){$11$}
\put(9,-1.5){$21$}
\put(12,-1.5){$2$}
\put(15,-1.5){$1$}
\put(18,-1.5){$\emptyset$}
\end{picture}
\end{center}
\end{example}
\vspace{1em}

\textsc{Remark.}  The \emph{Brauer algebra} $\mathfrak{B}_m$ is a
certain semisimple algebra, say over $\cc$, whose dimension is the
number $1\cdot 3\cdot\cdots(2m-1)$ of matchings on $[2m]$. (The
algebra $\mathfrak{B}_m$ depends on a parameter $x$ which is
irrelevant here.) Oscillating tableaux of length $2m$ are related to 
irreducible representations of $\mathfrak{B}_m$ in the same way that SYT of
content $[m]$ are related to irreducible representations of the
symmetric group $\sm$ and that vacillating tableaux of length $2m$ are
related to irreducible representations of the partition algebra
$\mathfrak{P}_m$. In particular, the irreducible representations
$J_\lambda$ of $\mathfrak{B}_m$ are indexed by partitions $\lambda$
for which there exists an oscillating tableau of shape $\lambda$ and
length $2m$, and $\dim J_\lambda$ is the number of such oscillating tableaux.
See e.g.\ \cite[Appendix~B6]{b-r} for further information.

Next we use the bijection between complete matchings and oscillating
tableaux to study the enumeration of $k$-noncrossing matchings. All
the results in the following hold for $k$-nonnesting matchings as
well.

For complete matchings, Theorem~\ref{lattice_P} becomes the following.
\begin{cor} \label{lattice}
The number of $k$-noncrossing matchings of $[2m]$ is equal to the
number of closed lattice walks of length $2m$ in the set
\[
V_k=\{ (a_1, a_2, \dots, a_{k-1}): a_1 \geq a_2 \geq \cdots \geq
a_{k-1}  \geq 0, a_i \in\zz \}
\]
from the origin to itself with unit steps in any coordinate direction
or its negative. ~~~~$\Box$ 
\end{cor}

Restricted to the cases $k=2, 3$, Corollary~\ref{lattice} leads to some
nice combinatorial correspondences.  Recall that a \emph{Dyck path} of
length $2m$ is a lattice path in the plane from the origin $(0,0)$ to
$(2m, 0)$ with steps $(1,1)$ and $(1, -1)$, that never passes below
the $x$-axis. A pair $(P, Q)$ of Dyck paths is \emph{noncrossing} if
they have the same origin and the same destination, and $P$ never goes
below $Q$.

\begin{cor} \label{lattice23}
\begin{enumerate}
\item The set of 2-noncrossing matchings is in one-to-one
  correspondence with  the set of Dyck paths.
\item The set of 3-noncrossing matchings is in one-to-one
correspondence with the set of  pairs of noncrossing Dyck paths.
\end{enumerate}
\end{cor}
\noindent {\bf Proof}.
By Corollary~\ref{lattice},  $2$-noncrossing matchings are in
one-to-one correspondence with closed lattice paths
$\{\vec v_i=(x_i)\}_{i=0}^{2m}$
with $x_0=x_{2m}=0$, $x_i \geq 0$ and $x_{i+1}-x_i=\pm 1$.
Given such a 1-dimensional lattice path, define a lattice path in the
plane by
letting $P=\{ (i, x_i) \st  i=0, 1, \dots, 2m\}$.
Then $P$ is a Dyck path, and this gives the desired correspondence.

For $k=3$, $3$-noncrossing matchings  are in
one-to-one correspondence with 2-dimensional lattice paths $\{
\vec v_i=(x_i, y_i) \}_{i=0}^{2m}$
with $(x_0, y_0)=(x_{2m}, y_{2m})=(0,0)$, $x_i \geq y_i \geq 0$ and
$(x_{i+1}, y_{i+1})-(x_i, y_i)=(\pm 1, 0)$ or $(0, \pm 1)$. Given
such a lattice path, define two lattice paths in the plane by
setting $P=\{(i, x_i+y_i)  \st  i=0, 1, \dots, 2m\}$ and
 $Q=\{(i, x_i-y_i)  \st  i=0, 1, \dots, 2m\}$. Then $(P, Q)$ is a pair of
noncrossing Dyck paths. It is easy to see that this is a bijection.
~~~$\Box$

\begin{example}
We illustrate the bijections between $3$-noncrossing matchings,
oscillating tableaux, and pairs of noncrossing Dyck paths.
The oscillating tableau is
\[
\eset, 1, 2, 21, 31, 21, 11,
21, 2, 1, \eset.
\]
 The sequence of SYT's as defined in the bijection is
\[
\emptyset, 1, 12,
\begin{array}{l} 1~2\\[-.05in] 3\end{array},
\begin{array}{l} 1~2~4\\[-.05in] 3\end{array},
\begin{array}{l} 1~2\\[-.05in] 3\end{array},
\begin{array}{l} 1\\[-.05in] 3\end{array},
\begin{array}{l} 1~7\\[-.05in] 3\end{array},\
37,\ 3,\ \emptyset.
\]
The corresponding matching and the pair of noncrossing
Dyck paths are given in
Figure~\ref{bijs}.
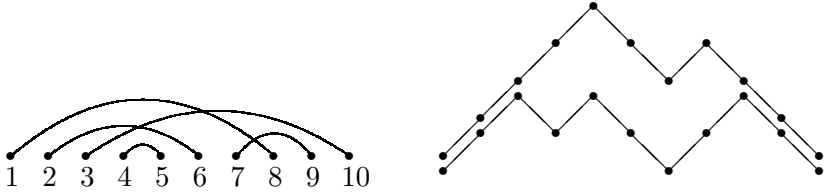
\begin{figure}[ht]
\begin{center}
\begin{picture}(380,80)
\setlength{\unitlength}{5mm}
\multiput(0.5,0)(1,0){10}{\circle*{0.2}}
\qbezier(0.5,0)(4,3)(7.5,0)\put(0.35,-0.8){1} \put(7.35,-0.8){8}
\qbezier(1.5,0)(3.5,1.6)(5.5,0)\put(1.35,-0.8){2}
\put(5.35,-0.8){6} \qbezier(2.5,0)(6,2.4)(9.5,0)\put(2.35,-0.8){3}
\put(9.3,-0.8){10} \qbezier(3.5,0)(4,0.6)(4.5,0)\put(3.35,-0.8){4}
\put(4.35,-0.8){5} \qbezier(6.5,0)(7.5,1.2)(8.5,0)
\put(6.35,-0.8){7}\put(8.35,-0.8){9}

\put(12,0){\line(1,1){1}}\put(12,0){\circle*{0.2}}
\put(13,1){\line(1,1){1}}\put(13,1){\circle*{0.2}}
\put(14,2){\line(1,1){1}}\put(14,2){\circle*{0.2}}
\put(15,3){\line(1,1){1}}\put(15,3){\circle*{0.2}}
\put(16,4){\line(1,-1){1}}\put(16,4){\circle*{0.2}}
\put(17,3){\line(1,-1){1}}\put(17,3){\circle*{0.2}}
\put(18,2){\line(1,1){1}}\put(18,2){\circle*{0.2}}
\put(19,3){\line(1,-1){1}}\put(19,3){\circle*{0.2}}
\put(20,2){\line(1,-1){1}}\put(20,2){\circle*{0.2}}
\put(21,1){\line(1,-1){1}}\put(21,1){\circle*{0.2}}
\put(22,0){\circle*{0.2}}
\put(12,-0.4){\line(1,1){1}}\put(12,-0.4){\circle*{0.2}}
\put(13,0.6){\line(1,1){1}}\put(13,0.6){\circle*{0.2}}
\put(14,1.6){\line(1,-1){1}}\put(14,1.6){\circle*{0.2}}
\put(15,0.6){\line(1,1){1}}\put(15,0.6){\circle*{0.2}}
\put(16,1.6){\line(1,-1){1}}\put(16,1.6){\circle*{0.2}}
\put(17,0.6){\line(1,-1){1}}\put(17,0.6){\circle*{0.2}}
\put(18,-0.4){\line(1,1){1}}\put(18,-0.4){\circle*{0.2}}
\put(19,0.6){\line(1,1){1}}\put(19,0.6){\circle*{0.2}}
\put(20,1.6){\line(1,-1){1}}\put(20,1.6){\circle*{0.2}}
\put(21,0.6){\line(1,-1){1}}\put(21,0.6){\circle*{0.2}}
\put(22,-0.4){\circle*{0.2}}
\end{picture}
\end{center}
\caption{The matching and the pair of noncrossing Dyck paths.}
\label{bijs}
\end{figure}
\end{example}

Let $f_k(m)$ be the number of $k$-noncrossing matchings of $[2m]$. By
Corollary~\ref{lattice} it is also the number of lattice paths of length
$2m$ in the region $V_k$ from the origin to itself with step set
$\{ \pm \epsilon_1, \pm \epsilon_2, \dots, \pm \epsilon_{k-1}\}$. Set
\[
F_k(x)=\sum_m f_k(m) \frac{x^{2m}}{(2m)!}.
\]
It turns out that a determinantal expression for $F_k(x)$ has been given
by Grabiner and Magyar \cite{GM93}.
It is simply the case $\lambda =
\eta = (m,m-1,...,1)$ of equation (38) in \cite{GM93}, giving
\begin{eqnarray}\label{gm}
 F_k(x) = \det\left[ I_{i-j}(2x) - I_{i+j}(2x) \right]_{i,j=1}^{k-1},
\end{eqnarray}
where
\[
   I_m(2x) = \sum_{j\geq 0}\frac{ x^{m+2j}}{j!(m+j)!},
\]
the hyperbolic Bessel function of the first kind of order $m$ \cite{WW27}.
One can easily check that when $k=2$,  the
generating function of $2$-noncrossing matchings equals
$$F_2(x)=I_0(2x)-I_2(2x)=\sum_{j\geq 0}C_j\frac{x^{2j}}{(2j)!},$$
where $C_j$ is the $j$-th Catalan number. When $k=3$, we have
\begin{eqnarray*}
f_3(m) &=&\frac{3!(2m+2)!}{m!(m+1)!(m+2)!(m+3)!}\\
&=&C_{m}C_{m+2}-C_{m+1}^2.
\end{eqnarray*}
This result agrees
with the formula on the number  of pairs of noncrossing Dyck paths
due to Gouyou-Beauchamps in \cite{GB89}.

{\sc Remark}.
The determinant formula \eqref{gm}  has been studied by  Baik and Rains
in \cite[Eqs. (2.25)]{BR00}. One simply puts $i-1$ for $j$ and
$j-1$ for $k$ in (2.25) of \cite{BR00} to get our formula. 
The same  formula was also obtained by Goulden \cite{Goulden} 
as the generating function for fixed point free permutations with no 
decreasing subsequence of length greater than $2k$. See Theorem 1.1 and 2.3 
of \cite{Goulden} and specialize $h_i$ to be $x^i/i!$, so $g_l$ becomes the 
hyperbolic Bessel function.  
The asymptotic distribution of $\mathrm{cr}(M)$ follows from another
result of Baik and Rains. In Theorem 3.1 of \cite{BR01} they obtained the
limit distribution for the length of the longest decreasing
subsequence of fixed point free involutions $w$. But representing $w$
as a matching $M$, the condition that $w$ has no decreasing
subsequence of length $2k+1$ is equivalent to the condition that $M$
has no $k+1$-nesting, and we already know that $\mathrm{cr}(M)$ and 
$\mathrm{ne}(M)$ have the same distribution. Combining the above
results, one has
\[  \lim_{m\rightarrow \infty} \mathrm{Pr}\left(
  \frac{\mathrm{cr}(M)-\sqrt{2m}}{(2m)^{1/6}}\leq \frac{x}{2}\right) = F_1(x),
\]
where
  $$ F_1(x) = \sqrt{F(x)}\exp\left( \frac 12\int_x^\infty
          u(s)ds\right), $$
where $F(x)$ is the Tracy-Widom distribution and $u(x)$ the Painlev\'e
II function.

Similarly one can try to enumerate complete matchings of $[2m]$ with
no $(k+1)$-crossing and no $(j+1)$-nesting. By the oscillating tableau
bijection this is just the number of walks of length $2m$ from
$\hat{0}$ to $\hat{0}$ in the Hasse diagram of the poset $L(k,j)$,
where $\hat{0}$ denotes the unique bottom element (the empty
partition) of $L(k,j)$, the lattice of integer partitions whose shape
fits in a $k \times j$ rectangle, ordered by inclusion. Let
$g_{k,j}(m)$ be this number, and $G_{k,j}(x)=\sum_m g_{k,j}(m) x^{2m}$
be the generating function.

For $j=1$, the number $g_{k,1}(m)$ counts lattice paths
from $(0,0)$ to $(2m,0)$ with steps $(1, 1)$ or $(1, -1)$ that stay
between the lines $y=0$ and $y=k$. The evaluation of  $g_{k,1}(m)$ was
first considered by Tak\'acs in \cite{Takacs62} by a probabilistic argument.
Explicit formula and generating function for this case are well-known. 
For example, in \cite{Mohanty79} one obtains the explicit
formula by applying the reflection principle repeatedly, viz.,
\begin{eqnarray*}
g_{k,1}(m)=\sum_i \left[ \binom{2m}{m-i(k+2)}
  -\binom{2m}{m+i(k+2)+k+1}\right].
\end{eqnarray*}
The generating function $G_{k,1}(x)$ is a special case of the one for 
the duration of the game
in the classical ruin problem, that is, restricted random walks with
absorbing barriers at $0$ and $a$, and initial position $z$.
See, for example, Equation (4.11) of Chapter 14 of \cite{feller68}: Let
\[
U_z(x)=\sum_{m=0}^\infty u_{z,m}x^m,
\]
where $u_{z,n}$ is the probability that the process ends with the $n$-th step
at the barrier $0$. Then
\[
U_z(x)=\left(\frac{q}{p}\right)^z
\frac{\lam_1^{a-z}(x)-\lam_2^{a-z}(x)}{\lam_1^{a}(x)-\lam_2^{a}(x)},
\]
where
\[
\lam_1(x)=\frac{1+\sqrt{1-4pqx^2}}{2px}, \qquad
\lam_2(x)=\frac{1-\sqrt{1-4pqx^2}}{2px}.
\]
The generating function
 $G_{k,1}(x)$ is just $U_1(2x)/x$ with $a=k+2$, $z=1$, and $p=q=1/2$.

In general, by the transfer matrix method \cite[{\S}4.7]{ec1}
\begin{eqnarray}
G_{k,j}(x)=\frac{\det(I-xA_{k,j}(0))}{\det(I-xA_{k,j})}
\end{eqnarray}
is a rational function,  where $A_{k,j}$ is the adjacency matrix
of the Hasse diagram of $L(k,j)$, and $A_{k,j}(0)$ is  obtained
from $A_{k,j}$ by deleting the row and the column corresponding to
$\hat{0}$. $A_{k,j}(0)$ is also the adjacency matrix of the Hasse
diagram of $L(k,j)$ with its bottom element (the empty partition)
removed. Note that   $\det(I-xA_{k,j})$ is a polynomial in $x^2$
since $L(k,j)$ is bipartite \cite[Thm.~3.11]{c-d-s}. Let
$\det(I-xA_{k,j})=p_{k,j}(x^2)$. The following is a table of
$p_{k,j}(x)$ for the values of $1 \leq k \leq j \leq 4$. (We only
need to list those with $k \leq j$ since
 $p_{k,j}(x)=p_{j,k}(x)$.)
\begin{eqnarray*}
\begin{array}{|c|l|}
\hline
(k,j) & p_{k,j}(x)=\det(I-\sqrt{x}A_{k,j}) \\ \hline
(1,1) & 1-x \\
(1,2) & 1-2x \\
(1,3) & 1-3x+x^2 \\
(1,4) & (1-x)(1-3x) \\
(2,2) & (1-x)(1-5x)\\
(2,3) & (1-x)(1-3x)(1-8x+4x^2)) \\
(2,4) & (1-14x+   49x^2-49x^3)(1-6x+5x^2-x^3) \\
(3,3) & (1-x)(1-19x+83x^2-x^3)(1-5x+6x^2-x^3)^2\\
(3,4) &  (1-2x)^2(1-8x+8x^2)(1-4x+2x^2)^2(1-16x+60x^2-32x^3+4x^4) \\
      &   \hspace{.5cm} \cdot   (1-24x+136x^2-160x^3+16x^4)  \\
(4,4) &
    (1-x)^2(1-18x+81x^2-81x^3)^2(1-27x+99x^2-9x^3)(1-9x+18x^2-9x^3)^2 \\
     & \hspace{.5cm} \cdot
  (1-27x+195x^2-361x^3)(1-6x+9x^2-x^3)^2(1-9x+6x^2-x^3)^2   \\
\hline
\end{array}
\end{eqnarray*}
The polynomial  $p_{k,j}(x)$ seems to have a lot of factors. We are
grateful to Christian Krattenthaler for 
explaining equation
(\ref{eq:kratt}) below, from which we can explain the factorization of
$p_{k,j}(x)$. By an observation \cite[{\S}5]{grabiner} of Grabiner,
$g_{k,j}(n)$ is equal to the number of walks with $n$ steps $\pm e_i$
from $(j, j-1, \dots, 2, 1)$ to itself in the chamber
$j+k+1>x_1>x_2>\cdots>x_j>0$ of
the affine Weyl group $\tilde{C}_n$.
Write $m=j+k+1$. By
\cite[(23)]{grabiner} there follows
  $$ \sum_n  g_{k,j}(n)\frac{x^{2n}}{(2n)! } = \det\left[ \frac 1m \sum_{r=0}^{2m-1}
     \sin(\pi ra/m)\sin(\pi rb/m)
    \cdot \exp(2x\cos(\pi r/m))\right]_{a,b=1}^j.
  $$
When this determinant is expanded, we obtain a linear combination of
terms of the form
 \begin{equation} \exp(2x(\cos(\pi r_1/m)+\cdots+\cos(\pi r_j/m)))
     = \sum_{n\geq 0} 2^n(\cos(\pi r_1/m)+\cdots+\cos(\pi r_j/m))^n
      \frac{x^n}{n!}, \label{eq:exp} \end{equation}
where $0\leq r_i\leq 2m-1$ for $1\leq i\leq j$. 
In fact, the case $\eta=\lambda$ of Grabiner's formula
\cite[(23)]{grabiner} shows that the number of walks of length $n$ in
the Weyl chamber from \emph{any} integral point to itself is again a
linear combination of terms $2^n(\cos(\pi r_1/m)+\cdots+\cos(\pi
r_j/m))^n$.
It follows that every eigenvalue of $A_{k,j}$ has the
form 
 \begin{equation} \theta=2(\cos(\pi r_1/m)+\cdots+\cos(\pi
      r_j/m)). \label{eq:kratt} \end{equation}
(In particular, the Galois group over $\mathbb{Q}$ of every irreducible
factor of $p_{k,j}(x)$ is abelian.) Note that \emph{a priori} not
every such $\theta$ may be an eigenvalue, since it may appear with
coefficient $0$ after the linear combinations are taken. 
The algebraic integer $z = 2(\cos(\pi r_1/m) +\dots +\cos(\pi r_j /m))$ 
lies in the field $\mathbb{Q}(\cos(\pi/m))$, an extension of $\mathbb{Q}$
of degree $\phi(2m)/2$, where $\phi$ is the Euler phi-function. 
To see this, let $z$ be a primitive $2m$-th root of unity. Then $z$ is 
a root of $x + 1/x = 2 \cos(\pi/m)$.
Hence the field $L = \mathbb{Q}(z)$ is quadratic or linear over $K =
\mathbb{Q}(\cos(\pi/m)$. Since $K$ is real and $L$ is not for $m>1$, 
we cannot have $K=L$. Hence $[L:K] = 2$. 
Since $[L:\mathbb{Q}] = \phi(2m)$, we have $[K:\mathbb{Q}] = \phi(2m)/2$.
It follows that 
the minimal polynomial over $\mathbb{Q}$ of $z$ has degree dividing
$\phi(2m)/2$. Thus every irreducible factor of $\det(I-Ax)$ has degree
dividing $\phi(2m)/2$, 
explaining  why
 $p_{k,j}(x)$  has many factors. A more careful analysis
should yield more precise information about the factors of
$p_{k,j}(x)$, but we will not attempt such an analysis here.

An interesting special case of determining $p_{k,j}(x)$ is determining
its degree, since the number of eigenvalues of $A_{k,j}$ equal to 0 is
given by ${j+k\choose j}- 2\,\cdot\deg\, p_{k,j}(x)$. Equivalently, since
$A_{k,j}$ is a symmetric matrix, $2\cdot\deg\, p_{k,j}(x) =
\mathrm{rank}(A_{k,j})$. We have observed the following.

\begin{enumerate}
\item For  $k+j \leq 12$ and $1 \leq k \leq j$, $A_{k,j}$ is
invertible exactly
for $(k,j)=$$(1,1), (1,3), (1, 5), (1, 7)$, $(1, 9)$, $(1, 11)$, $(3,3)$,
$(3, 7)$, $(3, 9)$, $(5,5)$ and $(5, 7)$.
\item  $A_{1, j}$ is invertible if and only if $j$ is odd. This is true
because $L(1,j)$ is a path of length $j$, whose determinant satisfies
the recurrence $\det(A_{1,j})=-\det(A_{1,j-2})$. The statement follows
from the initial conditions $\det(A_{1,1})=-1$ and
$\det(A_{1,2})=0$.
\item If $A_{k,j}$ is invertible, then $kj$ is odd. To see this, 
let $X_0 (X_1)$ be the set of  integer partitions of even (odd)
$n$ whose shape fits in a $k \times j$ rectangle.
Since $L(k,j)$ is bipartite graph with vertex partition $(X_0, X_1)$,
a necessary condition for $A_{k,j}$ to be invertible is $|X_0|=|X_1|$.
That is, the generating function $\sum_n p(k,j,n)q^n=\mathbf{
k+j \choose k}$ must have a root at $q=-1$, 
where $p(k,j, n)$ is the number of integer partitions on $n$ whose shape
fits into a $k\times j$ rectangle. 
But the multiplicity of
$1+q$ in the Gaussian polynomial $\mathbf{k+j \choose k}$ is
$\lfloor \frac{k+j}{2} \rfloor -\lfloor \frac k2 \rfloor-\lfloor \frac j2
\rfloor$, which is $0$ unless both $j$ and $k$ are odd.

\end{enumerate}


Item 3 is also proved independently by Jason Burns, who also found a counterexample for the converse:  
For $k=3$ and $j=11$, $A_{3, 11}$ is not invertible, in fact its  corank is 6.
The invertibility of $A_{k,j}$ for $kj$ being odd is currently under 
investigation.

\vspace{2em} 

\noindent {\large \bf Acknowledgments} 

The authors would like to thank Professor
 Donald Knuth for carefully reading the 
manuscript and providing many helpful comments.

%
%

\end{document}